\documentclass[reqno]{amsart}

\usepackage{xcolor}
\newtheorem{Theorem}{Theorem}[section]

\newtheorem{Lemma}[Theorem]{Lemma}

\newtheorem{Remark}{Remark}[section]

\title{On extended  lifespan  for 1d damped wave equation}
\author[K. Fujiwara]{Kazumasa Fujiwara}

\address[K. Fujiwara]{
Graduate School of Mathematics, Nagoya University,
Furocho, Chikusaku, Nagoya,
464-8602, Japan
}

\email{fujiwara.kazumasa@math.nagoya-u.ac.jp}

\author[V. Georgiev]{Vladimir Georgiev}

\address[V.Georgiev]{
Department of Mathematics,
University of Pisa,
Largo Bruno Pontecorvo 5,
I - 56127 Pisa, Italy}
\address{
Faculty of Science and Engineering, Waseda University,
3-4-1, Okubo, Shinjuku-ku, Tokyo 169-8555, Japan}
\address{
Institute of Mathematics and Informatics,  Bulgarian Academy of Sciences, Acad. Georgi Bonchev Str., Block 8, Sofia, 1113, Bulgaria
}
\email{georgiev@dm.unipi.it}

\begin{document}


\begin{abstract}
In this manuscript,
a sharp lifespan estimate
of solutions to semilinear classical damped wave equation
is investigated
in one dimensional case
when the sum of initial position and speed is $0$ pointwisely.
Especially,
an extension of lifespan is shown in this case.
Moreover, existence of some global solutions are obtained by a direct computation.
\end{abstract}

\maketitle

\section{Introduction}
In this manuscript,
we study the Cauchy problem of the following classical damped wave equation:
	\begin{align}
	\begin{cases}
	\partial_t^2 u + \partial_t u - \Delta u = |u|^{p},
	&\ t \in (0,T_0), \quad x \in \mathbb R^n,\\
	u(0) = u_0,
	&\ x \in \mathbb R^n,\\
	\partial_t u(0) = u_1,
	&\ x \in \mathbb R^n.
	\end{cases}
    \label{eq:DW}
	\end{align}
In particular in this manuscript,
we are interested
in the lifespan estimates of solutions in one dimensional case.

The existence and nonexistence of time-global solutions to \eqref{eq:DW} has been extensively studied.
We define so-called Fujita exponent $p_F(n)$ by
	\[
	p_F(n) = 1 + \frac 2 n.
	\]
$p_F(n)$ is known to be the critical exponent to \eqref{eq:DW}, that is, criteria of the existence and nonexistence of time-global solutions with small data in the sense explained below.
Todorova and Yordanov \cite{TY01} showed that
for any $(f_0,f_1) \in W^{1,2} \times L^2$ supported compactly,
there exists $\varepsilon_0$ such that solutions to \eqref{eq:DW} exist globally with
$(u_0,u_1) = (\varepsilon f_0, \varepsilon f_1)$ for any $\varepsilon \in (0,\varepsilon_0)$
when
	\[
	p_F(n) < p
	\begin{cases}
	\leq \frac{n}{n-2} &\quad \mathrm{if} \quad n \geq 3,\\
	< \infty &\quad \mathrm{if} \quad n =1,2,
	\end{cases}
	\]
where solutions $u$ satisfy
    \[
    u \in C(\lbrack 0, \infty); W^{1,2}),\quad
    \partial_t u \in C(\lbrack 0, \infty); L^2).
    \]
Here $W^{m,p}$ means the usual Sobolev space defined by
    \[
    W^{m,p}
    = \{ f \in L^p \mid \partial^\alpha f \in L^p\ (|\alpha| \leq m)\}.
    \]
Moreover, it has been shown that
if $p \in (1, p_F(n))$
and if $(f_0, f_1) \in L^1 \times L^1$ is regular enough and satisfies
	\begin{align}
	\int_{\mathbb R^n} (f_0(x) + f_1(x)) dx  > 0
	\label{eq:2}
	\end{align}
then for sufficiently small $\varepsilon >0$,
the maximal existence time $T_0=T_0(u_0,u_1)$ of solutions with $(u_0,u_1) = (\varepsilon f_0, \varepsilon f_1)$ enjoys the estimates
	\[
	T_{n,p}(c_{f_0,f_1} \varepsilon)
	\leq T_0
	\leq T_{n,p}(C_{f_0,f_1}\varepsilon),
	\]
where
	\[
	T_{n,p}(\varepsilon)
	= \begin{cases}
	\varepsilon^{- \frac{2(p-1)}{2-n(p-1)}}
	& \mathrm{if} \quad p \in (1, p_F(n)),\\
	\exp(\varepsilon^{- \frac 2 n }) & \mathrm{if} \quad p = p_F(n)
	\end{cases}
	\]
and $c_{f_0,f_1}$ and $C_{f_0,f_1}$ are positive constants
depending on $(f_0,f_1)$ but independent of $\varepsilon$.
We list the related results in the following table:
\begin{table}[h]
\begin{tabular}{llll}
\hline
Li \& Zhou \cite{LZ95}
&$n \leq 2$
&$1 < p \leq p_F(n)$
&$T_0 < T_{n,p}(C \varepsilon)$\\
\hline
Nishihara \cite{N04}
&$n = 3$
&$1 < p \leq p_F(n)$
&$T_0 < T_{n,p}(C \varepsilon)$\\
\hline
Ikeda \& Wakasugi \cite{IW15}
&$n \geq 1$
&$1 < p < p_F(n)$
&$T_0 < T_{n,p}(C \varepsilon)$\\
\hline
Ikeda \& Ogawa \cite{IO16}
&$n \geq 1$
&$p = p_F(n)$
&$T_0 > T_{n,p}(c \varepsilon)$\\
\hline
Fujiwara, Ikeda, \& Wakasugi \cite{FIW19}
&$n \geq 1$
&$1 < p < p_F(n)$
&$T_0 > T_{n,p}(c \varepsilon)$\\
\hline
Ikeda \& Sobajima \cite{IS19}
&$n \geq 1$
&$p = p_F(n)$
&$T_0 < T_{n,p}(C \varepsilon)$\\
\hline
\end{tabular}
\end{table}

We also refer the reader to
\cite{IKTW20,KTW19,KQ02,Z01}
for related topic.

The condition \eqref{eq:2} is technically essential
to obtain the upper bound of $T_0$
by the analysis taken in \cite{IS19, IW15}.
Indeed, they focused on the weak form of \eqref{eq:DW}:
	\begin{align}
	&\int_0^T \int_{\mathbb R^n} u(t,x) (\partial_t^2 - \partial_t - \Delta) \phi(t,x) dx \ dt
	\nonumber\\
	&= \int_{\mathbb R^n} (u_0(x) + u_1(x)) \phi(0,x) dx
	+ \int_0^T \int_{\mathbb R^n} |u(t,x)|^p \phi(t,x) dx \ dt,
	\label{eq:4}
	\end{align}
where $\phi \in \mathcal S(\mathbb R \times \mathbb R^n; [0,\infty) )$ with $\partial_t \phi(0,x) = \phi(T,x) =0$ for any $x \in \mathbb R^n$.
Roughly speaking,
by choosing appropriate $\phi$, the LHS of \eqref{eq:4} is controlled by the second term on the RHS of \eqref{eq:4}
and this contradicts with the positivity of the first term of RHS
derived by \eqref{eq:2}.
On the other hand,
it is not generally known whether solutions exists globally or not
when $p \leq p_F(n)$ and the condition \eqref{eq:2} fails.

The purpose of this manuscript is to estimate lifespan for solutions to \eqref{eq:DW}
with non-trivial initial data $(u_0, u_1)$ satisfying $u_0 + u_1 = 0$,
where the condition \eqref{eq:2} fails.

 Next, we plan to explain why  new phenomena is expected in this case. When $(u_0,u_1)$ is small and regular enough
and satisfies $u_0+u_1 \neq 0$,
the corresponding free solutions approach asymptotically to $e^{t \Delta} (u_0+u_1)$,
where $e^{t\Delta}$ is heat semi-group.
On the other hand, if $u_0+u_1=0$, then free damped wave solutions decay faster
than free heat solutions in an appropriate sense.
See \cite[Theorem 1.6]{IIOW19} and references therein, for example.
Denoting by $\mathcal{A}$ the set of initial data $(f_0, f_1) \in L^1 \times L^1$ such that $f_0,f_1$ are regular enough and satisfy \eqref{eq:2}
and defining $T_{0}(\varepsilon f_0, \varepsilon f_1; p)$ as the lifespan of the solution with initial data $\varepsilon f_0, \varepsilon f_1$ with $(f_0,f_1) \in \mathcal{A}$ we can say that
	\begin{align}
	\lim_{\varepsilon \to 0} \frac{T_{0}( \varepsilon f_0, \varepsilon f_1 ; p)}{T_{0}( \varepsilon g_0, \varepsilon g_1;p)} = C >0
	\label{eq.mm1}
	\end{align}
for all $(f_0,f_1) \in \mathcal{A},$  $(g_0,g_1) \in \mathcal{A}$ and $1 < p < 3.$
For $p=3$ we have
 \begin{align}
	\lim_{\varepsilon \to 0} \frac{\log T_{0}( \varepsilon f_0, \varepsilon f_1 ; p)}{\log T_{0}( \varepsilon g_0, \varepsilon g_1 ; 3)} = C >0 .
	\label{eq.mm2}
	\end{align}
 The analysis of the time decay of the linear damped wave equation suggests that the assumption
 $$ u_0+u_1=0$$ might extend the  lifespan. To be more precise, we can define  $\mathcal{B}$ as the set of initial data $(f, -f), $ $ f \in L^1 \backslash \{0\},$ such that $f$ is regular enough. Then in the place of \eqref{eq.mm1} one should expect
\begin{align}
	 \frac{T_{0}( \varepsilon f, -\varepsilon f ; p)}{T_{0}( \varepsilon g_0, \varepsilon g_1;p)}\sim \varepsilon^{-R}, R>0, \varepsilon \to 0
	\label{eq.mm3}
	\end{align}
 provided $(f,-f) \in \mathcal{B},$  $(g_0,g_1) \in \mathcal{A}$ and $1 < p < 3.$ One might interpret \eqref{eq.mm3} as definition of the extended lifespan for solutions  with small data in the set $\mathcal{B}$ compared with the well - known life - span of solutions with initial data in $\mathcal{A}$. The exponent $R>0$ is a natural quantity characterizing this extended lifespan. In a similar way the extension power $R>0$ for the critical case $p=3$ can be found through the relation
\begin{align}
	\frac{\log T_{0}( \varepsilon f, -\varepsilon f ; 3)}{\log T_{0}( \varepsilon g_0, \varepsilon g_1 ; 3)} \sim \varepsilon^{-R}
	\label{eq.mm4}
	\end{align}
as $\varepsilon \to 0.$ In fact, such an extension effect is established in our first result.


\begin{Theorem}
\label{Theorem:T0}
Let $n=1$ and $1 < p \leq 3$.
For any $f \in (W^{1,1} \cap W^{1,p}) \backslash \{0\}$,
there exist $\varepsilon_f > 0$ and $c_f, C_f > 0$ such that
if $0 < \varepsilon < \varepsilon_f$,
$u_0 = \varepsilon f$,
and $u_1 = - \varepsilon f$,
then the corresponding mild solution
    \[
    u
    \in C([0,T_0); W^{1,1} \cap W^{1,p})
    \cap C^1((0,T_0); L^{1} \cap L^{p})
    \]
exists uniquely with
    \begin{align}
    T_{1,p}(c_f \varepsilon^p)
    \leq T_0
    \leq T_{1,p}(C_f \varepsilon^p).
    \label{eq:T0}
    \end{align}
\end{Theorem}
\begin{Remark}
   The above result shows that   the extended lifespan  relation \eqref{eq.mm3} is fulfilled with
    \begin{align}
	 R= \frac{2(p-1)^2}{3-p}
	\label{eq.mm6}
	\end{align}
 in the subcritical case $1 < p < 3.$ In the critical case $p=3$ we have \eqref{eq.mm4} with
 $$ R =2(p-1). $$
\end{Remark}

Mild solutions are defined to be a $C([0,T_0); L^1 \cap L^p)$ function satisfying
    \begin{align}
    u(t,x)
    &= S(t) (u_0+u_1)(x)
    + \partial_t S(t) u_0(x)
    + \int_0^t S(t-\tau)|u(\tau)|^p(x) d \tau.
    \label{eq:Ieq}
    \end{align}
When $n=1$, $S(t)$ is defined by
    \[
    S(t) f = \frac 1 2 \int_{-t}^t e^{-t/2} I_0 ( \frac{\sqrt{t^2-y^2}}{2} ) f(x-y) dy
    \]
and $I_0$ is the modified Bessel function of $0$ order.
For $\ell \geq 0$, $I_\ell$ may be defined by the following Taylor series:
    \[
    I_\ell(z) = \sum_{k=0}^\infty \frac{1}{k!(k+\ell)!} \bigg( \frac z 2 \bigg)^{2k+\ell}
    \tag*{\cite[8.445]{TOC}}
    \]

Theorem \ref{Theorem:T0}
may follow from the approximation
    \begin{align}
    u(t) \sim \partial_t S(t) u_0
    \label{eq:usim}
    \end{align}
holding for $t \in [0,T_1]$
when $u_0+u_1=0$
under a suitable condition,
where $T_1 = C_f \varepsilon^{1-p}$ and $u_0 = \varepsilon f$.
\eqref{eq:usim} implies
that the following two estimates hold:
    \begin{align*}
    \| u(T_1) \|_{W^{1,1} \cap W^{1,p}}
    + \| \partial_t u(T_1) \|_{L^1 \cap L^p}
    &\leq C \varepsilon^p,\\
    \int u(T_1,x) + \partial_t u(T_1,x) dx
    &\geq C \varepsilon^p.
    \end{align*}
We note that if $u_0+u_1 \neq 0$,
then a similar approximations
    \[
    u(t) \sim S(t) (u_0+u_1) \sim e^{t \Delta}(u_0+u_1)
    \]
may hold for $t \in [0, T_{1,p}(C \varepsilon_1)]$
under a suitable condition,
where
    \[
    \varepsilon_1 = \| u_0 \|_{W^{1,1} \cap W^{1,p}} + \| u_1 \|_{L^p}.
    \]
Therefore, the solution $u$ to \eqref{eq:DW} is expected to be approximated by
    \[
    u(t) \sim S(t-T_1) ( u(T_1) + \partial_t u(T_1))
    \]
for $t \in [T_1, T_1 + T_{1,p}(C \varepsilon^p)]$.
This is a sketch of the proof of Theorem \ref{Theorem:T0}.

The natural question arising from Theorems \ref{Theorem:T0}
may be
whether solutions exist globally or blow up with
initial data $(u_0,u_1)$ satisfying $u_0+u_1 < 0$.
In \cite[Theorem 1.2]{LZ95},
it is shown that
if $p>1$, $u_0 \equiv 0$, and $u_1 < 0$ is smooth and integrable,
then the corresponding solutions exist globally
and take non-positive value at any $t$ and $x$.
We note that \cite[Theorem 1.2]{LZ95} may be generalized in the following way
and gives a partial answer:
\begin{Theorem}
\label{Theorem:global}
Let $n=1$ and $p > 1$.
Let $u_0 \in W^{1.1} \cap W^{1,p}$,
and $u_1 \in (L^1 \cap L^p)$.
There exist $\varepsilon_1 > 0$ such that
if
    \[
    \| u_0 \|_{W^{1,1} \cap W^{1,p}}
    + \| u_1 \|_{L^{1} \cap L^p}
    \leq \varepsilon_1
    \]
and
    \[
    u_0(x),\ u_1(x) + \frac 1 2 u_0(x) \leq 0
    \]
for any $x \in \mathbb R$,
then there exists a non-positive valued mild solution
    \[
    u \in C([0,\infty); W^{1,1} \cap W^{1,p})
    \cap C^1((0,\infty); L^{1} \cap L^{p})
    \]
to \eqref{eq:DW}.
\end{Theorem}

In this manuscript,
we show the approximation \eqref{eq:usim}
and lifespan estimates in Section \ref{section:2}.
We remark that
in multi-dimensional case,
a similar approach seems to work.
However, the approximation may be shown
with a Fourier analysis of \cite{IIOW19}
under more complicated situation.
In Section \ref{section:3},
we give a proof of Theorem \ref{Theorem:global} without smooth approximation.

\section{Proof for Theorem \ref{Theorem:T0}}
\label{section:2}

\subsection{Proof of existence and lifespan estimate from below}
We first recall the estimate for kernel parts of $S(t)$ and its derivatives.
For simplicity, we denote
    \[
    \omega = \sqrt{t^2-y^2}
    \]
as long as there is no confusion.
We also denote $\langle x \rangle = \sqrt{1+x^2}$.

We rewrite $S(t)$ and its derivatives by
    \begin{align}
    S(t) f(x)
    &= \widetilde{\mathcal K}_0(t) \ast f (x),
    \label{eq:S}\\
    \partial_x S(t) f(x)
    &= e^{-t/2} \frac{f(x+t)-f(x-t)}{2} + \widetilde{\mathcal K}_1(t) \ast f(x),
    \label{eq:Sx}\\
    \partial_t S(t) f(x)
    &= e^{-t/2} \frac{f(x+t)+f(x-t)}{2}
    + \widetilde{\mathcal K}_2(t) \ast f(x),
    \label{eq:St}\\
    \partial_{t,x} S(t) f(x)
    &= e^{-t/2} \frac{f'(x+t)+f'(x-t)}{2}
    \nonumber\\
    & - e^{-t/2} \bigg( \frac{t}{16} + \frac 1 4 \bigg) (f(x+t)-f(x-t))
    + \widetilde{\mathcal K}_3(t) \ast f(x),
    \label{eq:Stx}\\
    \partial_t^2 S(t) f(x)
    &= e^{-t/2} \frac{f'(x+t)-f'(x-t)}{2}
    \nonumber\\
    &+ e^{-t/2} \bigg( \frac{t}{16} - \frac 1 2 \bigg)
    (f(x+t)+f(x-t))
    + \widetilde{\mathcal K}_4(t) \ast f(x),
    \label{eq:Stt}
    \end{align}
where
    \begin{align*}
    \mathcal K_0(t,y)&= \frac 1 2 e^{-t/2} I_0(\frac \omega 2),\\
    \mathcal K_1(t,y)
    &= - \frac{e^{-t/2}}{4} \frac{y}{\omega} I_1(\frac \omega 2),\\
    \mathcal K_2(t,y)&= \frac 1 4 e^{-t/2}
    \bigg( \frac{t}{\omega} I_1(\frac \omega 2)
    - I_0(\frac \omega 2) \bigg),\\
    \mathcal K_3(t,y)
    &= - e^{-t/2}
    \bigg( \frac{yt}{16 \omega^2} I_2( \frac \omega 2 )
    - \bigg( \frac{y}{8\omega} + \frac{yt}{4\omega^3} \bigg) I_1(\frac \omega 2)
    + \frac{yt}{16 \omega^2} I_0( \frac \omega 2 )
    \bigg),\\
    \mathcal K_4(t,y)&= e^{-t/2}
    \bigg(
    \frac{t^2}{16 \omega^2} I_2 (\frac \omega 2)
    - \bigg( \frac{t}{4 \omega} + \frac{y^2}{4 \omega ^3} \bigg) I_1(\frac \omega 2)
    + \bigg( \frac 1 8 + \frac{t^2}{16 \omega^2} \bigg) I_0 (\frac \omega 2)
    \bigg),
    \end{align*}
$I_1$ and $I_2$ are first kind modified Bessel functions of first and second order,
and
    \[
    \widetilde{\mathcal K_j}(t,y)
    = \begin{cases}
    \mathcal K_j(t,y) & \mathrm{if} \quad |y| \leq t,\\
    0 & \mathrm{if} \quad |y| \geq t
    \end{cases}
    \]
for $0 \leq j \leq 4$.
For $0 \leq j \leq 4$, $\mathcal K_j$ is estimated as follows.

\begin{Lemma}[{\cite[Theorem 1.2]{MN03}}]
\label{Lemma:Kernel}
For $|y| \leq t$ and $0 \leq j \leq 4$, we have
    \[
    |\mathcal K_j(t,y)| \leq C \frac{e^{-y^2/8t}}{\langle t \rangle^{(j+1)/2}}.
    \]
In particular,
for $1 \leq  q \leq \infty$ and $0 \leq j \leq 4$,
the following estimates hold:
    \[
    \| \widetilde{\mathcal K_j}(t) \ast h \|_{L^q}
    \leq C \langle t \rangle^{-j/2-1/2q'}
    \| h \|_{L^1}.
    \]
\end{Lemma}

Now we estimate $u$ with $(u_0,u_1) = \varepsilon (f, -f)$.
We rewrite \eqref{eq:Ieq} with \eqref{eq:S} and \eqref{eq:St} by
$u=\Phi_f(u)$, where $\Phi_f$ is given by
    \begin{align}
    \Phi_f (u)(t,x)
    &:= \varepsilon e^{-t/2} \frac{f(x+t) + f(x-t)}{2}
    + \varepsilon \widetilde{\mathcal K}_2(t) \ast f(x)
    \nonumber\\
    &+ \int_0^t \widetilde{\mathcal K}_0(t-\tau) \ast |u(\tau)|^p(x) d \tau.
    \label{eq:u}
    \end{align}
Differentiating \eqref{eq:u}
and rewriting the resulting equality
with \eqref{eq:Sx}, \eqref{eq:St}, \eqref{eq:Stx}, and \eqref{eq:Stt},
we obtain
    \begin{align*}
    \partial_x \Phi_f(u)(t,x)
    &= \varepsilon e^{-t/2} \frac{f'(x+t)+f'(x-t)}{2}\\
    &- \varepsilon e^{-t/2} \bigg( \frac{t}{16} + \frac 1 4 \bigg)
    (f(x+t)-f(x-t))\\
    &+ \varepsilon \widetilde{\mathcal K}_3(t) \ast f(x)\\
    &+\int_0^t e^{-(t-\tau)/2}
        \frac{|u(\tau,x+t-\tau)|^p-|u(\tau,x-t+\tau)|^p}{2} d \tau\\
    &+ \int_0^t \widetilde{\mathcal K}_1(t-\tau) \ast |u(\tau)|^p(x) d \tau.
    \end{align*}
and
    \begin{align*}
    \partial_t \Phi_f(u)(t,x)
    &= \varepsilon e^{-t/2} \frac{f'(x+t)-f'(x-t)}{2}\\
    &+ \varepsilon e^{-t/2} \bigg( \frac{t}{16} - \frac 1 2 \bigg)
    (f(x+t)+f(x-t))\\
    &+ \varepsilon \widetilde{\mathcal K}_4(t) \ast f(x)\\
    &+\int_0^t e^{-(t-\tau)/2}
        \frac{|u(\tau,x+t-\tau)|^p+|u(\tau,x-t+\tau)|^p}{2} d \tau\\
    &+ \int_0^t \widetilde{\mathcal K}_2(t-\tau) \ast |u(\tau)|^p(x) d \tau.
    \end{align*}
Set
    \[
    \| u \|_{X(t)}
    = \sup_{0 \leq \tau \leq t}
    \sum_{0 \leq k+\ell \leq 1} \sum_{q=1,p}
    \langle \tau \rangle^{1+\ell+k/2+1/2q'} \|\partial_t^\ell \partial_x^k u(\tau)\|_{L^q}.
    \]
Lemma \ref{Lemma:Kernel} implies that
the estimates
    \begin{align*}
    &\|\partial_t^\ell \partial_x^k \Phi_f(u)(t)\|_{L^q}\\
    &\leq C \varepsilon \langle t \rangle^{-1-\ell-k/2-1/2q'}
    + C \|u\|_{X(t)}^p
    \int_0^t \langle t-\tau \rangle^{-\ell-k/2-1/2q'}
    \langle \tau \rangle^{-(3p-1)/2} d\tau\\
    &\leq C \varepsilon \langle t \rangle^{-1-\ell-k/2-1/2q'}
    + C
    \langle t \rangle^{-\ell-k/2-1/2q'} \|u\|_{X(t)}^p
    \end{align*}
hold for $t \geq 0$, $\ell, k =0,1$ with $\ell + k \leq 1$, and $q=1, p$,
where we note that
    \[
    \frac{3p-1}{2} = p(1+\frac{1}{2p'}) > 1.
    \]
In order to obtain the second estimate above,
we have used \cite[Lemma3.1]{S68}.
For $j=1,2$,
we also have used the Gagliardo-Nirenberg inequality to obtain
    \[
    \langle t \rangle^{(3p-1)/2} \||u(t)|^p\|_{L^p}
    \leq \langle t \rangle^{(3p-1)/2} \|u(t)\|_{L^{p}} \|\partial_x u(t)\|_{L^1}^{p-1}
    \leq \| u \|_{X(t)}^p.
    \]
Therefore, the estimate
    \[
    \| \Phi_f(u) \|_{X(t)}
    \leq C \varepsilon + C \langle t \rangle \| u \|_{X(t)}^p
    \]
holds for $t \geq 0$.
This and similar estimate for $ \Phi_f(u) - \Phi_f(v) $ imply that
$\Phi_f$ is contraction map on
    \[
    \{ u \in C([0,T_1];\thinspace W^{1,1} \cap W^{1,p})
        \cap C^1((0,T_1);\thinspace L^{1} \cap L^{p}) \mid \|u\|_{X(T_1)}
        \leq C \varepsilon \}
    \]
with $T_1 = c_f \varepsilon^{1-p}$.
We note that
we have
    \begin{align}
    \|u(T_1)\|_{W^{1,1} \cap W^{1,p}}
    + \|\partial_t u(T_1)\|_{L^1 \cap L^p}
    \leq T_1^{-1} \|u\|_{X(T_1)}
    \leq C \varepsilon^p.
    \label{eq:EstimateT1}
    \end{align}
Next we construct solutions from $t=T_1$.
Put
    \[
    \| u \|_{Y(t)}
    = \sup_{T_1 \leq \tau \leq t}
    \sum_{0 \leq k+\ell \leq 1} \sum_{q=1,p}
    \langle \tau - T_1 \rangle^{\ell+k/2+1/2q'} \|\partial_t^\ell \partial_x^k u(\tau)\|_{L^q}.
    \]
We rewrite \eqref{eq:Ieq} for $t \geq T_1$ by $u = \Psi(u)$,
where
    \begin{align*}
    \Psi(u)(t)
    &= S(t-T_1) (u(T_1)+ \partial_t u(T_1))
    + \partial_t S(t-T_1) u(T_1)\\
    &+ \int_{T_1}^t S(t-\tau)|u(\tau)|^p d \tau.
    \end{align*}
Then by a similar computation,
\eqref{eq:EstimateT1} implies that the estimate
    \[
    \| \Psi (u) \|_{Y(t)}
    \leq C \varepsilon^p
    + \| \Psi (u) \|_{Y(t)}^p \int_{T_1}^t \langle \tau-T_1 \rangle^{-\frac{p-1}{2}} d \tau
    \]
holds for $t \geq T_1$.
Therefore, $\Psi$ is shown a contraction map on
    \[
        \{ u \in C([T_1,T_2];\thinspace W^{1,1} \cap W^{1,p})
        \cap C^1((T_1,T_2);\thinspace L^{1} \cap L^{p}) \mid \|u\|_{Y(T_2)}
        \leq C \varepsilon_1 \},
    \]
    where $\varepsilon_1 = c_f \varepsilon^p $
    \[
    T_2-T_1 \sim \varepsilon_1^{2(p-1)/(3-p)} \sim T_{1,p} (c_f \varepsilon^p).
    \]
This shows the first estimate of \eqref{eq:T0}.

\subsection{Proof of upper bound lifespan estimate}
Since $u$ is constructed in $X(T_1)$,
the estimate
    \[
    \sup_{0\leq t \leq T_1} \bigg\| \int_0^t  S(t-\tau) |u(\tau)|^p d \tau \bigg\|_{L^p}
    \leq C \varepsilon^p
    \]
holds.
Moreover, since $S(\cdot)f \in C([0,\infty);L^p)$,
there exists $T_f \in (0,T_1)$ such that
    \[
    \inf_{0 \leq t \leq T_f} \| \partial_t S(t) f \|_{L^p}
    \geq \frac 1 2 \| f \|_{L^p}.
    \]
Since $T_f$ is independent of $\varepsilon$,
\eqref{eq:Ieq} implies that we have
    \[
    \inf_{0 \leq t \leq T_f} \| u(t) \|_{L^p}
    \geq C \varepsilon.
    \]
Since a straightforward computation implies that
    \[
    \int S(t) h(x) dx
    = (1-e^{-t}) \int h(x) dx
    \]
holds,
integrating \eqref{eq:Ieq} and differentiated \eqref{eq:Ieq},
we have
    \[
    \int ( \partial_t u(T_1,x) + u(T_1,x) ) dx
    = \int_0^{T_1} \| u(\tau)\|_{L^p}^p  d \tau
    \geq C T_f \varepsilon^p.
    \]
Therefore, combining the arguments of \cite{IS19,IW15},
we have
    \[
    T_0 \leq T_f + T_{1,p}(C \varepsilon^p).
    \]

\section{Proof of Theorem \ref{Theorem:global}}
\label{section:3}
In this section,
we show an a priori estimate
    \begin{align}
    &S(t) (u_0+u_1) + \partial_t S(t) u_0
    \nonumber\\
    &\leq u(t)
    \leq
    e^{-t/2} \frac{u_0(x-t)+u_0(x+t)}{2}
    + \frac{e^{-t/2}}{2} \int_{x-t}^{x+t}
    \Big(u_1(y) + \frac 1 2 u_0(y) \Big) dy
    \label{eq:apiroi}
    \end{align}
for any $t \geq 0$.
\eqref{eq:apiroi} and \eqref{eq:u} imply that we have the following uniform control
    \begin{align*}
    &| \partial_x u(t,x) - \partial_x S(t) (u_0+u_1)(x)
    + \partial_{x}\partial_{t} S(t)(u_0)(x) |\\
    &+| \partial_t u(t,x) - \partial_t S(t) (u_0+u_1)(x)
    + \partial_{t}^2 S(t)(u_0)(x) |\\
    &\leq \int_0^t e^{-(t-\tau)/2} |S(t) (u_0+u_1)(x) + S(t)(u_0)(x)|^p d\tau\\
    &+ \int_0^t (|\widetilde{\mathcal K}_1|+ |\widetilde{\mathcal K}_2|)
        \ast |S(t) (u_0+u_1) + S(t)(u_0)|^p(x) d\tau.
    \end{align*}
Namely, it is shown by the argument of Section \ref{section:2} that
there exists $A \geq 1$ such that
the estimate
    \begin{align}
    \| u(t) \|_{Y(t)} \leq A \varepsilon_1
    \label{eq:Yapriori}
    \end{align}
holds as long as \eqref{eq:apiroi} holds.
The existence of global solutions
follows from this and a blowup-alternative argument.

In \cite{LZ95},
\eqref{eq:apiroi} is shown by rewriting \eqref{eq:Ieq} by
    \begin{align}
    u(t,x)
    &= e^{-t/2} \frac{u_0(x-t)+u_0(x+t)}{2}
    + \frac{e^{-t/2}}{2} \int_{x-t}^{x+t}
    \Big(u_1(y) + \frac 1 2 u_0(y) \Big) dy
    \nonumber\\
    &+ \frac 1 2 \int_0^t \int_{x-t+\tau}^{x+t-\tau} e^{(\tau-t)/2}
    \Big( \frac{1}{4} + |u(\tau,y)|^{p-2} u(\tau,y) \Big) u(\tau,y) dy \ d \tau
    \label{eq:uWaveForm}
    \end{align}
with classical solutions $u$.
We note that
\eqref{eq:uWaveForm} follows from the fact that
$v = e^{t/2} u$ satisfies the following wave equation:
    \[
    \partial_t^2 v - \partial_x^2 v
    = e^{t/2} ( \frac 1 4 u + |u|^p).
    \]

Here we show that
\eqref{eq:uWaveForm} holds for mild solutions
by a direct computation without smooth approximation.
We claim that
    \begin{align}
    S(t) h(x)
    = \frac 1 8 \int_0^t e^{(\tau-t)/2} \int_{x-t+\tau}^{x+t-\tau}
    S(\tau) h(\xi) d \xi d \tau
    + e^{-t/2} \frac 1 2 \int_{x-t}^{x+t} h(y) dy
    \label{eq:STrans}
    \end{align}
for any $h \in L^1$.
Then \eqref{eq:uWaveForm} follows
for $u \in C([0,T);L^1 \cap L^p)$
from a straightforward computation with \eqref{eq:STrans}.
Now we show \eqref{eq:STrans}.
We compute
    \begin{align*}
    &\int_0^t \int_{x-t+\tau}^{x+t-\tau} e^{\tau/2} S(\tau) h(\xi)
    d \xi \thinspace d \tau\\
    &= \frac 1 2 \int_0^t \int_{x-t+\tau}^{x+t-\tau} \int_{-\tau}^\tau
    I_0 ( \frac{\sqrt{\tau^2-\eta^2}}{2}) h(\xi - \eta)
    d \eta \thinspace d \xi \thinspace d \tau\\
    &= \frac 1 2 \int_0^t \int_{-t+\tau}^{t-\tau} \int_{-\tau}^\tau
    I_0 ( \frac{\sqrt{\tau^2-\eta^2}}{2}) h(x - \xi - \eta)
    d \eta \thinspace d \xi \thinspace d \tau.
    \end{align*}
We change integral variables by $(y,\alpha,\beta) = (\eta+\xi,\tau+\eta,\tau-\eta)$.
Then Jacobian is $1/2$ and we compute that
    \begin{align*}
    &\frac 1 2 \int_0^t \int_{-t+\tau}^{t-\tau} \int_{-\tau}^\tau
    I_0 ( \frac{\sqrt{\tau^2-\eta^2}}{2}) h(x - \xi - \eta)
    d \eta \thinspace d \xi \thinspace d \tau\\
    &=\frac 1 4 \int_{-t}^t \int_{0}^{t-y} \int_{0}^{t+y}
    I_0 ( \frac{\sqrt{\alpha \beta}}{2}) h(x - y)
    d \alpha \thinspace d \beta \thinspace d y\\
    &=\frac 1 4 \int_{-t}^t \int_{0}^{t-y} \int_{0}^{t+y}
    \sum_{k=0} \frac{2^{-4k}}{(k!)^2} \alpha^k \beta^k h(x - y)
    d \alpha \thinspace d \beta \thinspace d y\\
    &= 4 \int_{-t}^t
    \sum_{k=0} \frac{2^{-4(k+1)}}{(k+1!)^2} (t^2-y^2)^{k+1} h(x - y) dy\\
    &= 4 \int_{-t}^t
    \Big(I_0(\frac{\sqrt{t^2-y^2}}{2})-1\Big) h(x - y) dy.
    \end{align*}
This implies \eqref{eq:STrans}.
We remark that \eqref{eq:uWaveForm} implies that the estimates
    \begin{align*}
    &S(t) (u_0+u_1) + \partial_t S(t) u_0\\
    &\leq
    e^{-t/2} \frac{u_0(x-t)+u_0(x+t)}{2}
    + \frac{e^{-t/2}}{2} \int_{x-t}^{x+t}
    \Big(u_1(y) + \frac 1 2 u_0(y) \Big) dy
    \leq 0
    \end{align*}
hold under the assumption of Theorem \ref{Theorem:global}.

Now we show \eqref{eq:apiroi} by \eqref{eq:uWaveForm}.
If \eqref{eq:apiroi} holds with $\varepsilon_1$ small enough,
then \eqref{eq:uWaveForm} is rewritten by
    \begin{align}
    |u(t,x)|
    &= -e^{-t/2} \frac{u_0(x-t)+u_0(x+t)}{2}
    - \frac{e^{-t/2}}{2} \int_{x-t}^{x+t}
    \Big(u_1(y) + \frac 1 2 u_0(y) \Big) dy
    \nonumber\\
    &+ \frac 1 2 \int_0^t \int_{x-t+\tau}^{x+t-\tau} e^{(\tau-t)/2}
    \Big( \frac 1 4 |u(\tau,y)| - |u(\tau,y)|^{p} \Big) dy \ d \tau.
    \label{eq:absuWaveForm}
    \end{align}
By a standard contraction argument,
solution to \eqref{eq:absuWaveForm} is constructed in
    \[
    \{ u \in C([0,T_2] ; L^1 \cap L^\infty) \mid
    \sup_{t \in [0,T_2]} (\| u(t) \|_{L^1} + \langle t \rangle \| u(t) \|_{L^\infty}) \leq
    C A \varepsilon_1\}
    \]
with $T_2 = T_{1,p}(C A \varepsilon_1)$.
By taking $\varepsilon_1$ small enough,
the uniqueness of solutions implies
that $u=-|u|$, \eqref{eq:apiroi}, and \eqref{eq:Yapriori} hold for any $t \in [0,T_2]$.
Repeating this argument with time interval whose width is $T_2$,
\eqref{eq:apiroi} is shown for any $t \geq 0$.

\section*{Acknowledgment}
The first author was supported in part by
JSPS Grant-in-Aid for Early-Career Scientists No. 20K14337. The second author was partially supported by   Gruppo Nazionale per l'Analisi Matematica, by the project PRIN  2020XB3EFL with the Italian Ministry of Universities and Research, by Institute of Mathematics and Informatics, Bulgarian Academy of Sciences, by Top Global University Project, Waseda University and the Project PRA 2022 85 of University of Pisa.

\end{document}